\documentclass{article}

\usepackage[ansinew]{inputenc}
\usepackage{amsmath}
\usepackage{amssymb}
\usepackage{amsthm}
\usepackage{amscd}
\usepackage{amsfonts}
\usepackage[english]{babel}
\usepackage{graphicx}
\newtheorem {defn}{Definition}

\newtheorem {thm}[defn]{Theorem}

\newtheorem {lemma}[defn]{Lemma}

\newcommand{\tr}{\operatorname{tr}}
\newcommand{\bbmatrix}[1]{\left[ \begin{array}{cccc} #1 \end{array} \right]}

\newcommand{\SP}{{\operatorname{sp}}}
\newcommand{\TR}{{\operatorname{tr}}}
\newcommand{\SPEC}{{\operatorname{spec}}}
\newcommand{\diag}{\operatorname{diag}}
\newcommand{\Rstar}{\boldsymbol{\mathcal R}_*}

\title{On completely integrable polynomial PDEs arising from Sturm-Liouville differential equation
using evolutionary vessels. KdV Hierarchy}
\author{Andrey Melnikov\\
Drexel University, Philadelphia, USA}

\begin{document}

\maketitle
\abstract{In this work we present a scheme for construction of solutions for evolutionary PDEs of some polynomial types $q'_t = P(q,q'_x,\ldots)$, where $P$ is a polynomial 
in a finite number of variables. This scheme is a generalization of the existing technique for solution of completely integrable 
PDEs using Inverse Scattering of the Sturm-Liouville differential equation.
The KdV equation 
\[ q'_t = - \dfrac{3}{2} q q'_x + \dfrac{1}{4} q'''_{xxx}
\]
is a special case, corresponding to type 1 evolutionary equations.
We present a complete solution of type 0, and present a KdV hierarchy corresponding to infinite number of polynomial
evolutionary equations rather for $\beta = \dfrac{1}{2} \int_0^x q(y,t)dy$ then for $q(x,t)$ itself, of the form
\[ \beta'_t = i^n b_n(\beta_x'),
\]
where  $b_0 = -\dfrac{1}{4} \beta'''_{xxx} + \dfrac{3}{2} (\beta'_x)^2$ corresponds to the KdV equation 
and $4 (b_{n+1})'_x = -i (b_n)_{xxx}''' + 4i (\beta'_xb_n)'_x$. Soliton solutions (i.e. involving pure exponents only) are presented for each such evolutionary equation,
demonstrating a ``simplicity'' of the solutions construction.
}

\tableofcontents

\section{Introduction}
Solution of the KdV equation plays a very significant role in the study of evolutionary Partial Differential equations (PDEs). Recall that the equation of the form
\[ q'_t = f(q,q'_x,\ldots)
\]
is called an \textit{evolutionary PDE} for the function of two variables $q(x,t)$ and a ``nice'' function $f$. Here $q'_x=\dfrac{\partial}{\partial x} q(x,t)$ and
$q'_t=\dfrac{\partial}{\partial x} q(x,t)$. Some times the solution $q(x,t)$ is also called a \textit{flow}.
In particular case we obtain the Korteweg-de Vries equation
\begin{equation}\label{eq:KdV}
q'_t = - \dfrac{3}{2} q q'_x + \dfrac{1}{4} q'''_{xxx}.
\end{equation}
named after two mathematicians D. J. Korteweg and  G. de Vries \cite{bib:KdV}. 
There is enormous amount of constructive examples of solutions for the KdV equation:
using algebraic-geometry data \cite{bib:Krich77}, flows on manifolds \cite{bib:LaxPeriod, bib:LaxAlmPeriod}, flows on Grassmanians 
\cite{bib:Sato, bib:BKSWSolitons, bib:Hirota, bib:TanDate}, flows on Grassmanians using loop-groups 
\cite{bib:SegalWilson}.

Usually the differential equation \eqref{eq:KdV} is studied with $q(x,0)$ defined for all values of $x\in\mathbb R$, belonging to a particular class of functions. 
This is the so called the \textit{initial value problem} of
the KdV equation. A standard technique \cite{bib:BDT}, \cite{bib:GGKM} to solve this equation is using the Inverse Scattering theory \cite{bib:GL, bib:FadeevInv} 
of the Sturm Liouville (SL) differential equation with a spectral parameter $\lambda$:
\begin{equation} \label{eq:SL}
-y''_{xx} + q y = \lambda^2 y.
\end{equation}
In this technique, one associates to the parameter $q(x,0)$, called \textit{potential} a scattering data (or a spectral function $\rho(\lambda)$), which differs from class
to class. Evolving in time the scattering data using a simple differential equation (yielding $\rho(\lambda,t)$), one then
translates back to the corresponding to this data potential $q(x,t)$. It turns out that $q(x,t)$ satisfies the KdV equation \eqref{eq:KdV}
(see introduction of \cite{bib:BDT} and \cite{bib:GGKM} for the idea).

Using the idea of P. Lax of so called \textit{Lax pair} one can show that defining linear operators
\[ P = -\dfrac{\partial^2}{\partial x^2} + q(x,t), \quad B = \dfrac{\partial^3}{\partial x^3} + b_1 \dfrac{\partial}{\partial x}+b_2 \]
acting on the space of Schwartz functions, one finds that the Lax condition, for an arbitrary function $w$ in the Schwartz space
\[ (\dfrac{\partial}{\partial t} P - [P,B]) w = 0
\]
is equivalent to $b_1=-\dfrac{3}{2}q$, $b_2=-\dfrac{3}{4}q'_x$ and the evolution of the function $q(x,t)$ according to the KdV equation \eqref{eq:KdV}.
This explains the choice of the differential equation, we are dealing with from this point of view. On the other hand, the same equation arises when
one constructs solutions of the KdV equation using evolutionary vessels, which we further discuss.

Using a recently developed idea of a SL vessel $\mathfrak{V}$, the author was able \cite{bib:GenVessel} to implement the scattering technique for SL equation in a manner, which unifies 
many previous classical works \cite{bib:GL, bib:FadeevInv}. 
There was later discovered a notion of evolutionary vessel $\mathfrak{V}_{KdV}$ in \cite{bib:KdVVessels}, generalizing the SL vessel $\mathfrak{V}$,
solving the KdV equation \eqref{eq:KdV}.
The vessel $\mathfrak{V}_{KdV}$ is just a collection of operators and speces, where operators satisfy some differential and algebraic
equations, among which we could mention the Lyapunov equation \eqref{eq:Lyapunov}, and a variation of the wave equation \eqref{eq:nDBt}.

In this work we extend the original notion of the KdV vessel, generating a technique for solution of a series of evolutionary PDEs, defined by
\begin{equation} \label{eq:nKdV} \beta'_t = i^n b_n(\beta_x'), \quad \beta = \dfrac{1}{2} \int_0^x q(y,t)dy
\end{equation}
where  $b_0 = -\dfrac{1}{4} \beta'''_{xxx} + \dfrac{3}{2} (\beta'_x)^2$ correspond to the KdV equation 
and $4 (b_{n+1})'_x = -i (b_n)_{xxx}''' + 4i (\beta'_xb_n)'_x$. The simplicity of the construction of a solution for \eqref{eq:nKdV} is demonstrated by presenting
a formula of a soliton in Section \ref{sec:Solitons}.

Let us emphasize few points related to the presented results. The recursively defined differential equations \eqref{eq:nKdV} are part of the KdV hierarchy,
because such evolutions commute between them (Lemma \ref{lemma:nKdVCommute}). We present in this work a technique to study non commuting evolutions (Section
\ref{sec:GenEvol}) but this requires an additional research. Finally notice that it is easy to construct a vessel in $N\leq\infty$ variables so that
its $\beta(x,t_1,t_2,\ldots,t_N)$ will satisfy equations \eqref{eq:nKdV} for each variable.
For this to hold one has to require that the operators of the vessel will evolve with respect to the variable $t_n$ as the operators of the $n$-th KdV vessel
\eqref{eq:nDBt}, \eqref{eq:nDXt}. There is a hope that it will enable to reprove the Witten conjecture \cite{bib:Witten} using these new ideas.

\section{KdV Vessels}
In this section we present KdV vessels and their properties. Each vessel has four important notions, related to it:
the tau function, the transfer function, the moments and the differential ring $\mathfrak{\boldsymbol R_*}$. All of them are defined in the
next sections along with their basic properties.

\subsection{Definition}
We define a set of matrices, which will be frequently referred to. We notice that the theory presented here can be modified for other choice
of these matrices. For example, there exists a definition realizing Non Linear Schr\" odinger equation (see \cite[Conclusions and Remarks]{bib:KdVVessels}).
\begin{defn}\label{def:VesPar}
SL vessel parameters for the SL equation \eqref{eq:SL} are defined as follows
\[ \sigma_1 = \bbmatrix{0 & 1 \\ 1 & 0},
\sigma_2 = \bbmatrix{1 & 0 \\ 0 & 0},
\gamma = \bbmatrix{0 & 0 \\ 0 & i},
\]
\end{defn}
Using these parameters we define a KdV vessel.
\begin{defn} 
A KdV evolutionary vessel, associated to the SL vessel parameters is a collection of operators and spaces
\begin{equation} \label{eq:DefV}
\mathfrak{V}_{KdV} = (A, B(x,t), \mathbb X(x,t); \sigma_1, 
\sigma_2, \gamma, \gamma_*(x,t);
\mathcal{H},\mathbb C^2), 
\end{equation}
where the operators $A,\mathbb X(x,t):\mathcal H\rightarrow\mathcal H$, $B(x,t):\mathbb C^2\rightarrow\mathcal H$ satisfy the following vessel conditions (a.e.):
\begin{eqnarray}
\label{eq:DB} \frac{\partial}{\partial x} B  &  & = - (A \, B \sigma_2 + B \gamma)\sigma_1^{-1}, \\
\label{eq:DX} \frac{\partial}{\partial x} \mathbb X &  & =  B \sigma_2 B^*, \\
\label{eq:DBtKdV} & \dfrac{\partial}{\partial t} B & = i A \dfrac{\partial}{\partial x} B, \\ 
\label{eq:DXtKdV} & \dfrac{\partial}{\partial t} \mathbb X & = 
 i A B \sigma_2 B^* - i B \sigma_2 B^* A^* +
 i B\gamma B^*, \\
\label{eq:Linkage}
&  \gamma_*  & =  \gamma + \sigma_2 B^* \mathbb X^{-1} B \sigma_1 -
 \sigma_1 B^* \mathbb X^{-1} B \sigma_2, \\
\label{eq:Lyapunov} & A \, \mathbb X + \mathbb X\,  A^*  & = - B \sigma_1 B^*, \\
\label{eq:X=X*} & \mathbb X^* & = \mathbb X.
\end{eqnarray}
The operator $A$ may be unbounded, but is assumed to generate a $C_0$ semi-group.
\end{defn}
Details on the properties of the vessel $\mathfrak{V}_{KdV}$ may be found in \cite{bib:KdVVessels}. Notice that there is a standard technique for constructing
of such a vessel. Suppose that we are given a $2\times 2$ function of complex variable $\lambda$, realized \cite{bib:bgr} in the following form
\[ S(\lambda) = I - B_0^*\mathbb X_0^{-1} (\lambda I - A)^{-1}B_0\sigma_1, \text{ $I$ - identity $2\times 2$ matrix}
\]
for which the operators $B_0:\mathbb C^2$, $\mathbb X_0, A: \mathcal H\rightarrow\mathcal H$ satisfy
\[ \mathbb X_0=\mathbb X_0^*, \quad A \mathbb X_0 + \mathbb X_0 A^* + B_0\sigma_1B_0^* = 0.
\]
Then follow these steps:
\begin{enumerate}
	\item Define $B(x)$ as the solution of \eqref{eq:DB} with the initial value $B_0$,
	\item Define $B(x,t)$ as the solution of \eqref{eq:DBtKdV} with the initial value $B(x)$ just constructed. 
	Notice that for a diagonal operator $A$, which is very common in these types of realizations,
	the equation \eqref{eq:DBtKdV} is a wave equation,
	\item Define $\mathbb X(x) = \mathbb X_0 + \int_0^x B(y)\sigma_2B^*(y)dy$,
	\item Solve \eqref{eq:DXtKdV} with the initial value $\mathbb X(x)$.
\end{enumerate}
So, the vessel conditions \eqref{eq:DB},\eqref{eq:DX}, \eqref{eq:DBtKdV} and \eqref{eq:DXtKdV} are satisfied by the construction. The equations
\eqref{eq:Lyapunov} and \eqref{eq:X=X*} are so called \textit{permanent}, namely if they hold at $t=0, x=0$ then they hold for all $x,t$ and this is true
from the choice of $B_0, \mathbb X_0$. It worth noticing that the formula for $\mathbb X(x,t)$ can be explicitly written, since the second mixed derivatives, applied to $\mathbb X(x,t)$
are equal, which can be easily checked:
\begin{thm} Let $B(x,t)$ be a solution of the partial differential equations \eqref{eq:DB}, \eqref{eq:DBtKdV}.
Then the formula for $\mathbb X(x,t)$ is as follows
\begin{multline} \label{eq:Xformula}
\mathbb X(x,t) = \mathbb X_0 + \int_0^x B(y,t)\sigma_2B^*(x,t) dy + \\
+ \int_0^t [i A B(0,s) \sigma_2 B^*(0,s) - i B(0,s)\sigma_2 B^*(0,s) A^* +
 i B(0,s)\gamma B^*(0,s)] ds.
\end{multline}
\end{thm}
As a result, it is a simple check that the formula for $\mathbb X(x,t)$ satisfies \eqref{eq:DX}, \eqref{eq:DXtKdV}. The main reason why we consider and call
the vessel $\mathfrak{V}_{KdV}$ as a KdV vessel is the following
\begin{thm}\label{thm:KdVequation}
Let $\mathfrak{V}_{KdV}$ be a KdV vessel, then the potential 
\[ q(x,t) =-2\dfrac{\partial^2}{\partial x^2}\ln(\det\mathbb X_0^{-1}\mathbb X(x,t))\]
satisfies KdV equation \eqref{eq:KdV}:
\[ q'_t = - \dfrac{3}{2} q q'_x + \dfrac{1}{4} q'''_{xxx}.
\]
\end{thm}
See \cite{bib:KdVVessels} for the proof, but the main idea behind it is that using the equations of the vessel $\mathfrak{V}_{KdV}$ one can differentiate
$\int_0^x q(y)dy =-2\TR(\sigma_2B^*(x,t)\mathbb X^{-1}(x,t)B(x,t)$ (see later in the text equation \eqref{eq:beta}) to verify the KdV equation \eqref{eq:KdV}.

\subsection{The transfer function and the tau-function}
By definition, the \textit{transfer function} of the vessel $\mathfrak{V}_{KdV}$ is defined as follows:
\begin{equation}\label{eq:DefS} 
S(\lambda,x,t) = I - B^*(x,t) \mathbb X^{-1}(x,t) (\lambda I - A)^{-1} B(x,t) \sigma_1,
\end{equation}
Notice that poles and singularities of $S$ with respect to $\lambda$ are determined by $A$ only.
Multiplication by the function $S(\lambda,x)$ maps \cite{bib:MelVin1}, \cite{bib:PerQuasiPerVessels} solutions of the input LDE with the spectral parameter $\lambda$
\begin{equation}
\label{eq:InCC}
 \lambda \sigma_2 u(\lambda, x) -
\sigma_1 \frac{\partial}{\partial x}u(\lambda,x) +
\gamma u(\lambda,x) = 0.
\end{equation}
to solution of the output LDE with the same spectral parameter
\begin{equation} \label{eq:OutCC}
\lambda \sigma_2 y(\lambda, x) - \sigma_1(x) \frac{\partial}{\partial x}y(\lambda,x) +
\gamma_* y(\lambda,x) = 0.
\end{equation}
The function $\gamma_*$ is defined by the Linkage condition \eqref{eq:Linkage}.
So, it follows that the function $S(\lambda,x)$ realizes a B\" acklund transformation of the SL equation (see \cite{bib:PerQuasiPerVessels}
for details).
As a result the transfer function satisfies the following ODE:
\begin{equation}\label{eq:DS} 
\dfrac{\partial}{\partial x} S(\lambda,x) =
\sigma_1^{-1}(\sigma_2\lambda+\gamma_*(x)) S(\lambda,x) -
S(\lambda,x) \sigma_1^{-1}(\sigma_2\lambda+\gamma).
\end{equation}
From the Lyapunov equation \eqref{eq:Lyapunov}, using simple calculations one can show that
\begin{eqnarray}
\label{eq:SSym}  S^*(-\bar\lambda,x,t)\sigma_1 S(\lambda,x,t)=\sigma_1, \\
\label{eq:SSym*} S(\lambda,x,t) \sigma_1^{-1} S^*(-\bar\lambda,x,t) = \sigma_1^{-1}.
\end{eqnarray}
Following the ideas presented in \cite{bib:SLVessels} we define the tau function of the transfer function
$S(\lambda,x)$ in the following way
\begin{defn} \label{def:Tau} For a given realization 
\[ S(\lambda,x,t) = I - B^*(x,t) \mathbb X^{-1}(x,t) (\lambda I - A)^{-1} B(x,t) \sigma_1
\]
tau function $\tau(x)$ is defined as
\begin{equation} \label{eq:Tau} \tau = \det (\mathbb X_0^{-1} \mathbb X(x,t))
\end{equation}
\end{defn}
Notice that from formula \eqref{eq:Xformula} it follows that
\begin{multline*}
\mathbb X_0^{-1} \mathbb X(x,t) = I + \mathbb X^{-1}(0) \int_0^x B(y,t)\sigma_2B^*(x,t) dy + \\
+ \mathbb X_0^{-1} \int_0^t [i A B(0,s) \sigma_2 B^*(0,s) - i B(0,s)\sigma_2 B^*(0,s) A^* +
 i B(0,s)\gamma B^*(0,s)] ds.
\end{multline*}
Since $\sigma_2$ has rank 1, this expression is of the form $I + T$, for a trace class operator $T$ and since 
$\mathbb X_0$ is an invertible operator, there exists a non trivial interval (of length at least $\dfrac{1}{\|\mathbb X_0^{-1}\|}$) on which $\mathbb X(x)$ and $\tau(x)$ are defined. Recall \cite{bib:GKintro} that a function $F(x)$ from (a, b) into the group G (the set of bounded invertible operators on H of the form I + T, for
a trace-class operator $T$) is said to be differentiable if $F(x) -I$ is \textit{differentiable} as a map into the trace-class operators. In our case,
\[ \dfrac{d}{dx} (\mathbb X_0^{-1}\mathbb X(x)) = 
\mathbb X_0^{-1} \dfrac{d}{dx} \mathbb X(x) =
\mathbb X_0^{-1} B(x)\sigma_2B^*(x)
\]
exists in trace-class norm. Israel Gohberg and Mark Krein \cite[formula 1.14 on p. 163]{bib:GKintro}
proved that if $\mathbb X_0^{-1}\mathbb X(x)$ is a differentiable function
into G, then $\tau(x) = \SP (\mathbb X_0^{-1}\mathbb X(x))$
\footnote{$\SP$ - stands for the trace in the infinite dimensional space.} is a differentiable map into $\mathbb C^*$ with
\begin{multline} \label{eq:GKform}
\dfrac{\tau'}{\tau}  = \SP (\big(\mathbb X_0^{-1}\mathbb X(x)\big)^{-1} 
\dfrac{d}{dx} \big(\mathbb X_0^{-1}\mathbb X(x)\big)) = \SP (\mathbb X(x)' \mathbb X^{-1}(x)) = \\
= \SP (B(x)\sigma_2 B^*(x) \mathbb X^{-1}(x)) =
\TR (\sigma_2 B^*(x) \mathbb X^{-1}(x)B(x)).
\end{multline}

Since any two realizations of a symmetric function
are weakly isomorphic, one obtains from standard theorems \cite{bib:bgr} 
in realization theory of analytic at infinity functions that they will
have the same tau function, i.e. this notion is independent of the realization we choose for the given function
$S(\lambda,x)$.
The following theorem was established for the finite dimensional case in \cite{bib:SLVessels}. We present
now its generalization for a \textit{regular} case, i.e. when all the operators of the vessel are bounded.
\begin{thm} \label{prop:Gamma*Formula}
For a Sturm Liouville regular vessel the following formula for $\gamma_*$ holds
\begin{equation} \label{eq:gamma*}
\gamma_* = \gamma + \bbmatrix{i\frac{\tau''}{\tau} & \frac{\tau'}{\tau} \\ -\frac{\tau'}{\tau} & 0}
\end{equation}
\end{thm}
\noindent\textbf{Proof:} If we denote $-\beta=\dfrac{\tau'}{\tau}$, then $\gamma_* = \bbmatrix{i \dfrac{\tau''}{\tau} & \dfrac{\tau'}{\tau} \\ -\dfrac{\tau'}{\tau} & i}$ 
then it is necessary to prove that $-\beta = \dfrac{\tau'}{\tau}$. Consider now the formula
\eqref{eq:GKform}
\begin{equation} \label{eq:TauBeta}
 \dfrac{\tau'}{\tau}  = \TR (\sigma_2 B^*(x,t) \mathbb X^{-1}(x,t)B(x,t)).
\end{equation}
Notice that the expression $B^*(x,t) \mathbb X^{-1}(x,t)B(x,t)$ is the first moment in the expansion of 
\[ S(\lambda,x,t) = I - \sum_{i=0}^\infty\dfrac{H_i\sigma_1}{\lambda^{i+1}}. \]
From the general formula of the first moment appearing in \cite[section 5]{bib:SchurVessels} (in this article the moments are
different from the moments in this article by multiplication on $\sigma_1^{-1}$ from the right)
\begin{equation}\label{eq:H0}
 H_0(x,t) = \bbmatrix{-\beta & -i\dfrac{\beta'-\beta^2}{2} \\ i \dfrac{\beta'-\beta^2}{2} & h_0^{22} }
\end{equation}
it follows that
\begin{equation} \label{eq:BetaTau}
 \dfrac{\tau'}{\tau} = \TR (\sigma_2 H_0(x)) = -\beta
\end{equation}
as desired.
\qed

As a result, from  \eqref{eq:GKform} and \eqref{eq:BetaTau} we obtain that
\begin{equation} \label{eq:beta}
\beta = -\tr(\sigma_2B^*(x) \mathbb X^{-1}(x)B(x))
\end{equation}

\subsection{The output differential ring $\mathfrak{\boldsymbol R_*}$ and Moment relations}
Consider now the moments $H_n(x,t)$ of the transfer function $S(\lambda,x,t)$, defined from the following formula
\begin{equation} \label{eq:DefHn}
S(\lambda,x,t) = 
I - B^*\mathbb X^{-1}(\lambda I - A)^{-1}B\sigma_1 = I - \sum_{n=0}^\infty \dfrac{ H_n(x,t)}{\lambda^{n+1}}\sigma_1,
\end{equation}
or more precisely
\begin{defn}\label{defn:Hn} The \textbf{$n$-th moment} $H_n(x,t)$ of a regular vessel $\mathfrak{V}$ is:
\[ H_n(x,t) = B^*(x,t) \mathbb X^{-1}(x,t) A^n B(x,t).
\]
\end{defn}
Moments will play a crucial role in the future research and we will need some of their basic properties.
\begin{thm} The following relations between the moments of a KdV vessel $\mathfrak{V}_{KdV}$ hold
\begin{eqnarray}
H_0^*-H_0=0,\\
\label{eq:MomSymmetry} H_n + (-1)^{n+1} H_n^* + \sum\limits_{i=0}^{n-1} (-1)^i H_i^* \sigma_1 H_{n-1-i}  = 0, \\
\label{eq:MomSymmetry*} H_n + (-1)^{n+1} H_n^* + \sum\limits_{i=0}^{n-1} (-1)^{n-1+i} H_i \sigma_1 H^*_{n-1-i}  = 0, \\
\label{eq:DHn} 	(H_n)'_x = \sigma_1^{-1}\sigma_2H_{n+1} - H_{n+1} \sigma_2\sigma_1^{-1} + \sigma_1^{-1} \gamma_* H_n - H_n \gamma\sigma_1^{-1}.
\end{eqnarray}
\end{thm}
\noindent\textbf{Proof: 1.} the first equality and \eqref{eq:MomSymmetry} follow from the symmetry condition \eqref{eq:SSym}. The equation
\eqref{eq:MomSymmetry*} follows from \eqref{eq:SSym*}. The last \eqref{eq:DHn} follows from the differential equation \eqref{eq:DS}. \qed

Plugging in the formula for $(H_n)'_x$ the expressions appearing in Definition \ref{def:VesPar} for the vessel parameters, and using \eqref{eq:gamma*}
rewritten in terms of $\beta$, we obtain that
\begin{equation} \label{eq:DH}
	(H_n)'_x = \bbmatrix{0&0\\1&0} H_{n+1} - H_{n+1} \bbmatrix{0&1\\0&0} + \bbmatrix{\beta&i\\-i(\beta'-\beta^2)&-\beta} H_n - H_n \bbmatrix{0&0\\i&0}.
\end{equation}
Remember that our goal is to construct (solutions of) evolutionary equation so that the derivative of $\beta'_t$ is a polynomial of derivatives of
$\beta$ with respect to $x$. So it is necessary to consider the following differential ring of functions:
\begin{defn} The differential ring $\boldsymbol{\mathcal R}_*$ is generated by $\beta$ and its derivatives with respect to $x$.
\end{defn}
The role of the ring $\boldsymbol{\mathcal R}_*$ in the simplest case ($\dim\mathcal H<\infty$) was investigated in \cite{bib:SLVessels}.

\subsection{KdV equation in matrix form}
Let us consider in more details the case of KdV vessel. Using vessel conditions one can check that the following formula holds:
\[ \begin{array}{lllllll}
\dfrac{\partial}{\partial t} H_n & = \dfrac{\partial}{\partial t}[B^*\mathbb X^{-1}A^n B] = (B^*)'_t\mathbb X^{-1}A^nB - B^* \mathbb X^{-1} \mathbb X'_t \mathbb X^{-1}A^nB+
B^*\mathbb X^{-1}A^n B'_t = \\
& = \text{using evolutionary conditions \eqref{eq:DBtKdV}, \eqref{eq:DXtKdV} } \\
& = (B^*)'_x(-iA^*) \mathbb X^{-1}A^nB - B^* \mathbb X^{-1} (iA\mathbb X'_x-i\mathbb X'_x A^* + iB\gamma B^*) \mathbb X^{-1}A^nB+
B^*\mathbb X^{-1}A^n (iA) B'_x = \\
& = \text{using \eqref{eq:DB}, \eqref{eq:Linkage} and \eqref{eq:Lyapunov}} \\
& = i \dfrac{\partial}{\partial x} H_{n+1} + i \dfrac{\partial}{\partial x} [H_0] \sigma_1 H_n,
\end{array} \]
Consequently,
\begin{equation}\label{eq:DHntKdV}
\dfrac{\partial}{\partial t} H_n = i \dfrac{\partial}{\partial x} H_{n+1} + i \dfrac{\partial}{\partial x} [H_0] \sigma_1 H_n.
\end{equation}
As a result we obtain that the following PDE holds
\begin{equation} \label{eq:PDEforS}
 \dfrac{\partial}{\partial t} S(\lambda,x,t) = i\lambda \dfrac{\partial}{\partial x} S(\lambda,x,t) +
i \dfrac{\partial}{\partial x} [H_0] \sigma_1 S(\lambda,x,t).
\end{equation}
Notice that $\dfrac{\partial}{\partial x} [H_0]$ has entries, which depend on $\beta(x,t)$ and its $x$-derivatives only:
\[ \dfrac{\partial}{\partial x} [H_0] = \bbmatrix{-\beta'_x & - \dfrac{i(\beta'_{xx}-2\beta\beta'_x)}{2} \\ \dfrac{ i(\beta''_{xx}+2\beta\beta'_x)}{2} & \dfrac{d}{dx} h_0^{21} }
\]
and we will see later that $\dfrac{d}{dx} h_0^{21}$ is in $\boldsymbol{\mathcal R}_*$.

The KdV equation \eqref{eq:KdV} arises also as the compatibility condition for the equality of mixed derivatives, applied to $S(\lambda,x,t)$.
More precisely from
\[ \dfrac{\partial}{\partial x} [\dfrac{\partial}{\partial t} S(\lambda,x,t)] = 
\dfrac{\partial}{\partial t} [\dfrac{\partial}{\partial x} S(\lambda,x,t)] 
\]
using \eqref{eq:DS}, \eqref{eq:PDEforS}, \eqref{eq:DefHn} and considering the constants (not depending on $\lambda$ coefficients) one obtains that the following formula must hold
\begin{equation} \label{eq:Dgamma*tKdV}
(\gamma_*)'_t = - i \gamma_* (H_0)'_x\sigma_1 + i \sigma_1 (H_0)''_{xx} \sigma_1 +i \sigma_1 (H_0)'_x \gamma_*.
\end{equation} 
When we consider this equation for the entries of $\gamma_*=\bbmatrix{-i(\beta'_x-\beta^2)&-\beta\\\beta&i}$ this becomes the regular KdV equation \eqref{eq:KdV}. For this it is
enough to focus on the $1,2$ entry.

\section{Different types of evolution}
In the next section we present a general scheme for the evolution of $\beta(x,t)$. This idea generates flows, which commute and do not commute with the KdV flow.
In Section \ref{sec:KdVHier} we present an infinite number of commuting with KdV flows.
\subsection{\label{sec:GenEvol}General Evolution}
Suppose that evolution of the operator $B(x,t)$ with respect to $t$ is determined by the following formula
\begin{equation} \label{eq:nDBt}
 \dfrac{\partial}{\partial t}(B\sigma_1) = A \sum\limits_{i=0}^n A^i B m_i = A (B m_0 + A B m_1 + \ldots + A^n B m_n), 
\end{equation}
where $m_i$ are $2\times 2$ matrices, which will be determined later. 
Let us find first a differential
equation for $\mathbb X$ so that the permanency of the Lyapunov equation \eqref{eq:Lyapunov} holds. Notice that when the condition
$\SPEC{A}\cap\SPEC{-A^*}=\emptyset$ is fulfilled the solution $\mathbb X(x,t)$ of \eqref{eq:Lyapunov} is unique and as a result the operator
$\mathbb X(x,t)$ is uniquely determined. On the other hand, when the condition $\SPEC{A}\cap\SPEC{-A^*}\neq\emptyset$ holds, we have to solve for $\mathbb X$
using \eqref{eq:DX} using initial condition, which does satisfy \eqref{eq:Lyapunov}. It turns out that an alternative approach may be applied and the formula
for $\mathbb X$ will be explicitly obtained. Differentiating Lyapunov equation \eqref{eq:Lyapunov} we find the following formula
\[ \begin{array}{llll}
A \mathbb X'_t + \mathbb X'_t A^* + B'_t \sigma_1 B^* + B \sigma_1 (B^*)'_t & =
A \mathbb X'_t + \mathbb X'_t A^* + A (\sum\limits_{i=0}^n A^i B m_i B^*)+\sum\limits_{i=0}^nB\bar m_i B^* (A^*)^iA^* \\
& = A(\mathbb X'_t + \sum\limits_{i=0}^n A^i B m_i B^*) + (\mathbb X'_t + \sum\limits_{i=0}^n B \bar m_i B^*(A^*)^i)A^*=0.
\end{array} \]
So, the derivative $\mathbb X'_t$ must be chosen in such a manner that the last equality holds for an arbitrary $A$ and $B$, satisfying \eqref{eq:DBtKdV}.
Consider special cases for different values of $i$:
\begin{enumerate}
	\item[0:] in this case, $B'_t = A Bm_0$ and $A(\mathbb X'_t+Bm_0B^*) + (\mathbb X'_t+B\bar m_0B^*)A^*=0$. A simple solutions of this is
		\[ \mathbb X'_t = - Bm_0B^* = -\mathbb Y_0, \quad m_0 = \bar m_0 . \]
	\item[1:] Starting from $B'_t = A (B m_0 + A B m_1)$ we obtain that 
	\[ A(\mathbb X'_t+Bm_0B^* + ABm_1B^*) + (\mathbb X'_t+B\bar m_0B^*+B\bar m_1 A^*B^*)A^* = 0
	\]
	We have to choose, denoting by $\mathbb Y_1$ the new term:
	\[ \mathbb X'_t = -\mathbb Y_0 - A B m_1 B^* - B \bar m_1 B^* A^* = -\mathbb Y_0 - \mathbb Y_1,
		\quad \bar m_0=m_0, \quad \bar m_1 = -m_1
	\]
	\item[2:] $B'_t = A (B m_0 + A B m_1 + A^2 B m_2)$. Then in order to satisfy
	\[ A(\mathbb X'_t+Bm_0B^* + ABm_1B^*+ A^2 B m_2B^*) + (\mathbb X'_t+B\bar m_0B^* + 
		B\bar m_1 A^*B^*+ B \bar m_2B^*(A^*)^2)A^* = 0,
	\]
	we have to require
	\[ \mathbb X'_t = -\mathbb Y_0 - \mathbb Y_1-\mathbb Y_2,
	\]
	where $\mathbb Y_2 = A^2Bm_2B^* - ABm_2B^*A^* + B\bar m_2B^* (A^*)^2$ for a self-adjoint $2\times 2$ matrix
	$m_2=\bar m_2$.
\end{enumerate}
Generalizing these examples, we obtain the following theorem.
\begin{thm} Let $\dfrac{\partial}{\partial t}(B\sigma_1) = A \sum\limits_{n=0}^n A^n B m_n$ \eqref{eq:DBtKdV}. 
Then in order to obtain the permanency of the Lyapunov equation \eqref{eq:Lyapunov}, it is sufficient to require $\bar m_n = (-1)^n m_n$,
\begin{equation}
\label{eq:defY} \mathbb Y_n = \sum\limits_{i=0}^n (-1)^i A^{n-i} B m_n B^* (A^*)^i,
\end{equation}
and
\begin{equation} \label{eq:nDXt}
\mathbb X'_t = - \sum\limits_{i=0}^n \mathbb Y_i.
\end{equation}
\end{thm}
\noindent\textbf{Proof:} Simple calculations show that $\mathbb Y_n^* = \mathbb Y_n$ and
\[ A \mathbb Y_n + \mathbb Y_n A^* + A^{n+1} B m_n B^* + B \bar m_n B^* (A^*)^n = 0.
\]
As a result,
\[ A \mathbb X'_t + \mathbb X'_t A^* + \dfrac{\partial}{\partial t} [B \sigma_1 B] = 0.
\]
Thus the permanency condition of the Lyapunov  equation holds. \qed

\begin{defn} An  $n$-th KdV vessel is the colloection
\[ \mathfrak{V}_{nKdV} = (A, B(x,t), \mathbb X(x,t); \sigma_1, 
\sigma_2, \gamma, \gamma_*(x,t);
\mathcal{H},\mathbb C^2), 
\]
where the operators act as for the ordinary KdV vessel and satisfy the equations \eqref{eq:DB}, \eqref{eq:DX}, \eqref{eq:Linkage}, \eqref{eq:Lyapunov}, \eqref{eq:X=X*}
and the two evolutionary equations \eqref{eq:nDBt} and \eqref{eq:nDXt}.
The number $n$ appearing in the definition of $\mathfrak{V}_{nKdV}$ is called the \textbf{evolutionary type}.
\end{defn}

Let us closer investigate the formula for the evolution of a moment $H_n$
\[ \begin{array}{lllllll}
\dfrac{\partial}{\partial t} H_N & = \dfrac{\partial}{\partial t}[B^*\mathbb X^{-1}A^N B] \\
& = \dfrac{\partial}{\partial t}[B^*] \mathbb X^{-1}A^N B + B^*\dfrac{\partial}{\partial t}[\mathbb X^{-1}]A^N B + B^*\mathbb X^{-1}A^N\dfrac{\partial}{\partial t}B \\
& = \sigma_1^{-1} \sum\limits_{i=0}^n \bar m_i B^* (A^*)^{i+1} \mathbb X^{-1} A^N B + B^* \mathbb X^{-1}(\sum\limits_{i=0}^n \mathbb Y_i)\mathbb X^{-1} A^N B +
B^*\mathbb X^{-1}A^N \sum\limits_{i=0}^n A^{i+1} B m_i \sigma_1^{-1} \\
& = \sigma_1^{-1} \sum\limits_{i=0}^n \bar m_i B^* (A^*)^{i+1} \mathbb X^{-1} A^N B + B^* \mathbb X^{-1}(\sum\limits_{i=0}^n \mathbb Y_i)\mathbb X^{-1} A^N B + 
\sum\limits_{i=0}^n H_{N+i+1} m_i \sigma_1^{-1}
\end{array} \]
Thus we obtain the following equation:
\begin{equation} \label{eq:dHNtGen}
\dfrac{\partial}{\partial t} H_N = \sigma_1^{-1} \sum\limits_{i=0}^n \bar m_i B^* (A^*)^{i+1} \mathbb X^{-1} A^N B + B^* \mathbb X^{-1}(\sum\limits_{i=0}^n \mathbb Y_i)\mathbb X^{-1} A^N B + 
\sum\limits_{i=0}^n H_{N+i+1} m_i \sigma_1^{-1}
\end{equation}
Examining the formulas for $\mathbb Y_i$ and using the Lyapunov equation \eqref{eq:Lyapunov}, it is possible to rewrite this formula in the following form, after multiplication by
$\sigma_1$ on the right:
\[ \dfrac{\partial}{\partial t} H_N \sigma_1= \sum\limits_{i=0}^n P_i H_{N+i}\sigma_1 + \sum\limits_{i=0}^n H_{N+i+1} \sigma_1 \sigma_1^{-1}m_i,
\]
where $P_i$ are non commutative polynomials of the first $n$ moments, their adjoints and the constant matrices $\sigma_1, \sigma_2, m_i, \bar m_i$.
Consequently one obtains that a formula for the evolution of $S(\lambda,x,t)$ is
\[\dfrac{\partial}{\partial t} S = (\sum\limits_{i=0}^n \lambda^i P_i) S - S \sum\limits_{i=0}^n \lambda^{i+1} \sigma_1^{-1} m_i + Q,
\]
where $Q$ is a polynomial of degree $n$ in $\lambda$ with coefficients as for $P_i$, which makes the equality to hold for positive powers of $\lambda$.

Since
\[  \begin{array}{lll} 
\dfrac{\partial}{\partial t} \dfrac{\partial}{\partial x} S & = \dfrac{\partial}{\partial t} [\sigma_1^{-1}(\sigma_2\lambda+\gamma_*)S-S\sigma_1^{-1}(\sigma_2\lambda+\gamma)] \\
& = \sigma_1^{-1}(\gamma_*)'_t S + \sigma_1^{-1}(\sigma_2\lambda+\gamma_*)S'_t-
S'_t  \sigma_1^{-1}(\sigma_2\lambda+\gamma) \\
& = \sigma_1^{-1}(\gamma_*)'_t S + \sigma_1^{-1}(\sigma_2\lambda+\gamma_*) [(\sum\limits_{i=0}^n \lambda^i P_i) S - S \sum\limits_{i=0}^n \lambda^{i+1} \sigma_1^{-1} m_i +Q]- \\
&  \hspace{2cm} - [(\sum\limits_{i=0}^n \lambda^i P_i) S - S \sum\limits_{i=0}^n \lambda^{i+1} \sigma_1^{-1} m_i + Q]\sigma_1^{-1}(\sigma_2\lambda+\gamma)
\end{array} \]
and
\[ \begin{array}{lll}
\dfrac{\partial}{\partial x} \dfrac{\partial}{\partial t} S 
& = \dfrac{\partial}{\partial x} [(\sum\limits_{i=0}^n \lambda^i P_i) S - S \sum\limits_{i=0}^n \lambda^{i+1} \sigma_1^{-1} m_i + Q]  \\
& = (\sum\limits_{i=0}^n \lambda^i (P_i)'_x) S + (\sum\limits_{i=0}^n \lambda^i P_i) S'_x - S'_x \sum\limits_{i=0}^n \lambda^{i+1} \sigma_1^{-1} m_i + Q'_x\\
& =  (\sum\limits_{i=0}^n \lambda^i (P_i)'_x) S + (\sum\limits_{i=0}^n \lambda^i P_i) [\sigma_1^{-1}(\sigma_2\lambda+\gamma_*)S-S\sigma_1^{-1}(\sigma_2\lambda+\gamma)] - \\
& \hspace{2cm} - [\sigma_1^{-1}(\sigma_2\lambda+\gamma_*)S-S\sigma_1^{-1}(\sigma_2\lambda+\gamma)] \sum\limits_{i=0}^n \lambda^{i+1} \sigma_1^{-1} m_i + Q'_x
\end{array} \]
From here, we will obtain that the equality of mixed derivatives:
\begin{equation} \label{eq:MixedDevEqual}
\dfrac{\partial}{\partial t} \dfrac{\partial}{\partial x} S = \dfrac{\partial}{\partial x} \dfrac{\partial}{\partial t} S,
\end{equation}
is equivalent to
\begin{multline} \label{eq:PartialEqualGen}
 [\sigma_1^{-1}(\gamma_*)'_t + \sigma_1^{-1}(\sigma_2\lambda+\gamma_*) (\sum\limits_{i=0}^n \lambda^i P_i) ] S +
 S [\sum\limits_{i=0}^n \lambda^{i+1} \sigma_1^{-1} m_i] \sigma_1^{-1}(\sigma_2\lambda+\gamma) 
 + \sigma_1^{-1}(\sigma_2\lambda+\gamma_*) Q - Q \sigma_1^{-1}(\sigma_2\lambda+\gamma) = \\
 = [\sum\limits_{i=0}^n \lambda^i (P_i)'_x + (\sum\limits_{i=0}^n \lambda^i P_i) \sigma_1^{-1}(\sigma_2\lambda+\gamma_*)] S +
 S \sigma_1^{-1}(\sigma_2\lambda+\gamma) [\sum\limits_{i=0}^n \lambda^{i+1} \sigma_1^{-1} m_i] +Q'_x.
\end{multline}
In the KdV case notice that $Q=0$ and there is a commutative relation
\[ [\sum\limits_{i=0}^n \lambda^{i+1} \sigma_1^{-1} m_i] \sigma_1^{-1}(\sigma_2\lambda+\gamma) =
\sigma_1^{-1}(\sigma_2\lambda+\gamma) [\sum\limits_{i=0}^n \lambda^{i+1} \sigma_1^{-1} m_i],
\]
since for the KdV case $\sum\limits_{i=0}^n \lambda^{i+1} \sigma_1^{-1} m_i = i \lambda \sigma_1^{-1}(\sigma_2\lambda+\gamma)$.
Then the condition \eqref{eq:MixedDevEqual}, multiplied by the inverse of $S$ on the right, is
\[ \sigma_1^{-1}(\gamma_*)'_t + \sigma_1^{-1}(\sigma_2\lambda+\gamma_*) (\sum\limits_{i=0}^n \lambda^i P_i) =
\sum\limits_{i=0}^n \lambda^i (P_i)'_x + (\sum\limits_{i=0}^n \lambda^i P_i) \sigma_1^{-1}(\sigma_2\lambda+\gamma_*).
\]
Here the coefficient of $\lambda^0$ (i.e. of constant) is the desired evolutionary equation
\[ \sigma_1^{-1}(\gamma_*)'_t = - \sigma_1^{-1} \gamma_* P_0 + (P_0)'_x + P_0\sigma_1^{-1} \gamma_*
\]
coinciding with \eqref{eq:Dgamma*tKdV} for $P_0=i\sigma_1(H_0)'_x\sigma_1$. Notice that we can rewrite the last equality as
\begin{equation} \label{eq:EvolKdVHirhy}
\sigma_1^{-1}(\gamma_*)'_t = (P_0)'_x + [P_0, \sigma_1^{-1} \gamma_*],
\end{equation}
where $[A,B]=AB-BA$ is the usual notation for the commutator of matrices.
In the more general case we obtain that the term $Q=\sum\limits_{i=0}^n \lambda^i Q_i$ is present and the evolutionary condition is more
complicated. For this we have to consider the coefficients of $\lambda^0$ at the equation \eqref{eq:PartialEqualGen}:
\[ \sigma_1^{-1}(\gamma_*)'_t + \sigma_1^{-1}\gamma_*P_0 + \sigma_1^{-1} \gamma_*Q_0 -Q_0\sigma_1^{-1}\gamma =
(P_0)'_x + P_0 \sigma_1^{-1}\gamma_* + (Q_0)'_x
\]
which is equivalent to
\begin{equation}\label{eq:EvolGenHirhy}
\sigma_1^{-1}(\gamma_*)'_t = (P_0+Q_0)'_x + [P_0,\sigma_1^{-1}\gamma_*] - \sigma_1^{-1} \gamma_*Q_0 + Q_0\sigma_1^{-1}\gamma. 
\end{equation}
The main disadvantage of these formulas is that it is difficult to obtain a polynomial evolutionary equation \eqref{eq:EvolGenHirhy}.
Usually the terms $P_0$ and $Q_0$ involves integrals. We will see in Section \ref{sec:KdVHier} that it can be overcome for a special choice of the 
matrices $m_i$.

\subsection{Evolutionary type zero}
In the first simple case, when $B'_t \sigma_1 = A (B m_0)$, we can differentiate with respect to $t$ 
equality \eqref{eq:beta} to obtain
\[ \begin{array}{llll}
 -\beta'_t& = \tr[\sigma_2 \sigma_1^{-1} \bar m_0 B^* A^*\mathbb X^{-1}B + m_0 \sigma_1^{-1} \sigma_2 B^* \mathbb X^{-1}AB + 
 						\sigma_2B^* \mathbb X^{-1}B m_0 B^* \mathbb X^{-1}B ] = \\
 					&	= \tr[\sigma_2 \sigma_1^{-1} \bar m_0 B^* (-\mathbb X^{-1}A-\mathbb X^{-1}B\sigma_1B^*\mathbb X^{-1})B + m_0 \sigma_1^{-1} \sigma_2 B^* \mathbb X^{-1}AB + \\
 					& \quad \quad \quad + \tr[\sigma_2B^* \mathbb X^{-1}B m_0 B^* \mathbb X^{-1} B] = \\
 					&	= \tr[(-\sigma_2 \sigma_1^{-1} \bar m_0 + m_0 \sigma_1^{-1} \sigma_2)B^* \mathbb X^{-1}A B -\sigma_2 \sigma_1^{-1} \bar m_0 
 									B^*\mathbb X^{-1}B\sigma_1B^*\mathbb X^{-1}B + \\
 					& \quad \quad \quad + \tr[\sigma_2B^* \mathbb X^{-1}B m_0 B^* \mathbb X^{-1} B] = \\
					& = \tr[\bbmatrix{m_{12}-\bar m_{12}&-m_{22}\\m_{22}&0} H_1]
					-\tr[\bbmatrix{\bar m_{12}& m_{22}\\0&0} H_0 \sigma_1 H_0] +\\
					& \quad \quad \quad + \tr[\sigma_2 H_0 m_0  H_0].
\end{array} \]
where $m_0 = \bar m_0 = \bbmatrix{m_{11}&m_{12}\\\bar m_{12}&m_{22}}$. Thus
\[ \begin{array}{llll}
-\beta'_t & =  (m_{12}-\bar m_{12}) \tr( \sigma_2 H_1) - m_{22} \tr(\bbmatrix{0&1\\-1&0} H_1) -
						\bbmatrix{\bar m_{12}&\bar  m_{22}} H_0 \sigma_1 H_0 \bbmatrix{1\\0}) + \\
					& \quad \quad \quad + \tr[\sigma_2B^* \mathbb X^{-1}B m_0 B^* \mathbb X^{-1} B] = \\
					& = (m_{12}-\bar m_{12})\tr( \sigma_2 H_1) - m_{22} \tr(\bbmatrix{0&1\\-1&0} H_1) + \beta m_{22} \tr(\bbmatrix{0&0\\0&1} H_0) + \\
					& \quad \quad \quad  + \bbmatrix{\bar m_{12}&0} H_0 \bbmatrix{i\dfrac{\tau''}{2\tau}\\\dfrac{\tau'}{\tau}}) + 
						\bbmatrix{0 & m_{22}} H_0 \bbmatrix{i\dfrac{\tau''}{2\tau}\\0} + \tr[\sigma_2H_0 m_0 H_0].
\end{array} \]
Since the last line in this expression is in $\boldsymbol{\mathcal R}_*$, we have to require that
\begin{equation} \label{eq:Type0Dbetat}
(m_{12}-\bar m_{12})\tr( \sigma_2 H_1) - m_{22} \tr(\bbmatrix{0&1\\-1&0} H_1) + \beta m_{22} \tr(\bbmatrix{0&0\\0&1} H_0) \in \boldsymbol{\mathcal R}_*.
\end{equation}
Using \eqref{eq:Lyapunov} it is immediate that the first equality in the next list holds. Then differentiating these equalities and using \eqref{eq:DH} we obtain other
lines:
\begin{eqnarray}
\tr(\sigma_1 H_0)  = - \tr(A + A^*) = const \in \mathfrak{\boldsymbol R_*}, \\
\tr(\sigma_2 H_0)  = - \beta \in \mathfrak{\boldsymbol R_*}, \\ 
\tr(\bbmatrix{0&1\\-1&0} H_0)  =  i\beta'_x - i\beta^2 \in \mathfrak{\boldsymbol R_*}, \\
\label{eq:H1_1} \tr(\sigma_2 H_1) - \tr(\bbmatrix{0&0\\0&i} H_0)  =
 \beta \tr(\bbmatrix{0&1\\0&0} H_0) - \dfrac{1}{2i}(\beta''_{xx}-4\beta'_x\beta+2\beta^3)\in \mathfrak{\boldsymbol R_*}, \\
\label{eq:H1_1-1} \tr(\bbmatrix{2i\beta&-1\\1&0}H_1) + \beta \tr(\bbmatrix{0&0\\0&1}H_0)  = \\
\hspace{1cm} = i(\beta'-\beta^2)\tr(\bbmatrix{0&1\\0&0}H_0) +\dfrac{-1}{4}(\beta'''_{xxx} - 6 (\beta')^2)\in \mathfrak{\boldsymbol R_*}
\end{eqnarray}
Taking into account formulas \eqref{eq:H1_1}, \eqref{eq:H1_1-1} we have to require in \eqref{eq:Type0Dbetat} that $m_{12}-\bar m_{12} = m_{22}=0$ in order to obtain
in \eqref{eq:Type0Dbetat} an expression belonging to $\mathfrak{\boldsymbol R_*}$ . Finally, the formula for $\beta'_t$ becomes (for real $m_{11}, m_{12}$):
\[ \begin{array}{llllll}
-\beta' & = \tr[\sigma_2H_0 m_0 H_0] = \bbmatrix{1&0} H_0 \bbmatrix{m_{11}& m_{12}\\m_{12}&0} H_0 \bbmatrix{1\\0} = \\
& = \bbmatrix{\dfrac{\tau'}{\tau}&-i\dfrac{\tau''}{2\tau}} \bbmatrix{m_{11}& m_{12}\\m_{12}&0} \bbmatrix{\dfrac{\tau'}{\tau} \\ i\dfrac{\tau''}{2\tau}} \\
& = m_{11} [\dfrac{\tau'}{\tau}]^2 - i \dfrac{\tau''}{2\tau} m_{12} \dfrac{\tau'}{\tau} + i \dfrac{\tau''}{2\tau} m_{12} \dfrac{\tau'}{\tau} \\
& = m_{11} [\dfrac{\tau'}{\tau}]^2 
\end{array} \]
in other words, we obtain the following theorem:
\begin{thm} Suppose that $B'_t = A B (m \sigma_2 + m_{12}\sigma_1)$ for real constants $m, m_{12}$. Then the evolutionary equation for
$\beta$ is:
\begin{equation} \label{eq:beta'tm0}
 \beta'_t = - m \beta^2.
\end{equation}
Moreover, $m_0=m \sigma_2 + m_{12}\sigma_1$ is the only possible equation, for which $B'_t = A B m_0$ and $\mathbb X'_t = - Bm_0B^*$ defines an evolutionary 
equation of the polynomial form.
\end{thm}
Notice that the solution of this equation is
\[ \beta(x,t) = \dfrac{1}{mt-\alpha(x)}
\]
for an arbitrary function $\alpha(x)$. For $t=0$ we obtain that $\beta(x,0)=\dfrac{1}{\alpha(x)}$ which must be equal for the beta function of a SL vessel. In other,
words, starting with a SL vessel, realizing $\beta(x)$, the corresponding evolutionary vessel of type $0$ will have beta equal to
\[ \beta(x,t) = \dfrac{1}{mt-\dfrac{1}{\beta(x)}} = \dfrac{\beta(x)}{m t \beta(x) - 1}.
\]

\subsection{Evolutionary KdV type 2}
Suppose that $B$ satisfies the following evolutionary equation
\begin{equation} \label{eq:BtKdV2}
 B'_t \sigma_1 = -A^2 B'_x \sigma_1 = A^3 B \sigma_2 + A^2 B \gamma,
\end{equation}
which mean that $m_2=\sigma_2, m_1=\gamma$. Then the following theorem holds
\begin{thm} The following formulas hold
\begin{eqnarray}
 (H_n)'_t = - (H_{n+2})'_x - (H_1)'_x \sigma_1 H_{n+1} - (H_0)'_x\sigma_1H_0\sigma_1 H_n, \\
 S'_t = -\lambda^2 S'_x - [\lambda (H_0)'_x\sigma_1 + (H_1)'_x \sigma_1 + (H_0)'_x\sigma_1H_0\sigma_1] S,
\end{eqnarray}
and $\gamma_*$ satisfies the evolutionary equation \eqref{eq:EvolKdVHirhy} with $P_0 =  (H_1)'_x \sigma_1 + (H_0)'_x\sigma_1H_0\sigma_1$.
\end{thm}
\noindent\textbf{Remark:} Notice that the expression $\sigma_1 P_0$ is indeed anti-self-adjoint (as $\gamma_*$ is such) because of \eqref{eq:MomSymmetry} for $n=1$.
Differentiating this expression we obtain
\[ (H_1)'_x + (H_1^*)'_x + (H_0)'_x \sigma_1 H_0 + H_0 \sigma_1 (H_0)'_x = 0,
\]
and from here multiplying by $\sigma_1$ on both sides we will obtain the desired anti-self-adjointness.

\noindent\textbf{Proof:} The formula or the moments is straightforward, by using the evolutionary condition \eqref{eq:BtKdV2} and derived from this evolutionary condition
\eqref{eq:DXtKdV}. From this formula it follows that the evolutionary equation of $S$ is
\[ S'_t = -\lambda^2 S'_x - [\lambda (H_0)'_x\sigma_1 + (H_1)'_x \sigma_1 + (H_0)'_x\sigma_1H_0\sigma_1] S + Q(\lambda).
\]
Comparing the coefficients of all non-negative powers of $\lambda$, we come to the conclusion that $Q(\lambda)=0$ and this proves the evolutionary formula for $S$.
Formula for the evolution of $\gamma_*$ is immediate from \eqref{eq:EvolKdVHirhy}.
\qed

Thus it remains to investigate the formula, in order to understand what kind of evolution is obtained in this manner. We will see later ($n=2$ of Theorem \ref{thm:nKdV}) 
that the evolutionary formula obtained in this case is of the form $\beta'_t=P(\beta'_x,\beta''_{xx},\ldots)$ for a polynomial $P$.

\section{\label{sec:KdVHier}KdV hierarchy}
\subsection{Definition and analysis of the Hierarchy}
Let us generalize the construction of Evolutionary KdV type 2. We will suppose that $B$ satisfies the following evolutionary equation
\begin{equation} \label{eq:BtKdVn}
 B'_t \sigma_1 = (\sqrt{-1}A)^n B'_x \sigma_1.
\end{equation}
In other words, since $B'_x \sigma_1= - AB\sigma_2 - B\gamma$ this requirement is equivalent to the definition
\[ m_n = -(\sqrt{-1})^n \sigma_2, \quad m_{n-1}= - (\sqrt{-1})^n \gamma.
\]
Notice first that if the operator $B$ corresponds to a multidimensional vessel, i.e. depends on $N$ (possibly infinite number of) variables
$B=B(x,t_1,t_2,\ldots,t_N)$, so that the evolution with respect to $t_i$ is defined by the formula \eqref{eq:BtKdVn} with the power
$n_i$, then
\[ (B'_{t_i})'_{t_j} = (((\sqrt{-1})A)^{n_i} B'_x)'_{t_j} = ((\sqrt{-1})A)^{n_i} ((\sqrt{-1})A)^{n_j} B''_{xx}
\]
and the same result holds for $(B'_{t_j})'_{t_i}$ since $((\sqrt{-1})A)^{n_i} ((\sqrt{-1})A)^{n_j} = ((\sqrt{-1})A)^{n_j} ((\sqrt{-1})A)^{n_i}$. Thus we obtain:
\begin{lemma}\label{lemma:nKdVCommute} Evolutions, defined by the formula \eqref{eq:BtKdVn} commute (including commutation with the derivative with respect to $x$).
\end{lemma}
\noindent\textbf{Proof:} we have seen that the commutation holds for the operator $B$. In a similar manner, it can be shown for $\mathbb X$,
for this we have to prove that (using \eqref{eq:nDBt}, \eqref{eq:defY} and $m_n = -(\sqrt{-1})^n \sigma_2, \quad m_{n-1}= - (\sqrt{-1})^n \gamma$)
\begin{multline*}
-(\mathbb X'_{t_1})'_{t_2} = \dfrac{\partial}{\partial t_2} [ \sum\limits_{i=0}^{n_1} (-1)^i A^{n_1-i} B (\sqrt{-1})^{n_1}\sigma_2 B^* (A^*)^i + 
\sum\limits_{i=0}^{n_1-1} (-1)^i A^{n_1-1-i} B (\sqrt{-1})^{n_1}\gamma B^* (A^*)^i] = \\
= \dfrac{\partial}{\partial t_1} [ \sum\limits_{i=0}^{n_2} (-1)^i A^{n_2-i} B (\sqrt{-1})^{n_2}\sigma_2 B^* (A^*)^i + 
\sum\limits_{i=0}^{n_2-1} (-1)^i A^{n_2-1-i} B (\sqrt{-1})^{n_2}\gamma B^* (A^*)^i] = -(\mathbb X'_{t_2})'_{t_1}.
\end{multline*}
Let us calculate the left hand side of this equality.
\[ \begin{array}{llll}
\dfrac{\partial}{\partial t_2} [ \sum\limits_{i=0}^{n_1} (-1)^i A^{n_1-i} B (\sqrt{-1})^{n_1}\sigma_2 B^* (A^*)^i + 
\sum\limits_{i=0}^{n_1-1} (-1)^i A^{n_1-1-i} B (\sqrt{-1})^{n_1}\gamma B^* (A^*)^i] = \\
= \sum\limits_{i=0}^{n_1} (-1)^i A^{n_1-i} [(\sqrt{-1}A)^{n_2} B'_x (\sqrt{-1})^{n_1}\sigma_2 B^* +  B (\sqrt{-1})^{n_1}\sigma_2 (B^*)'_x (-\sqrt{-1}A^*)^{n_2}]  (A^*)^i + \\
\quad + \sum\limits_{i=0}^{n_1-1} (-1)^i A^{n_1-1-i} [(\sqrt{-1}A)^{n_2} B'_x (\sqrt{-1})^{n_1}\gamma B^* +  B (\sqrt{-1})^{n_1}\gamma (B^*)'_x (-\sqrt{-1}A^*)^{n_2}] (A^*)^i 
\end{array} \]
Dividing this expression by $(\sqrt{-1})^{n_1+n_2}$ and using formula \eqref{eq:nDBt} we can obtain after regrouping, that this expression is
\[ \begin{array}{llll}
\sum\limits_{i=0}^{n_1} (-1)^i A^{n_1-i} [A^{n_2} B'_x \sigma_2 B^* +  B \sigma_2 (B^*)'_x (-A^*)^{n_2}]  (A^*)^i + \\
\quad + \sum\limits_{i=0}^{n_1-1} (-1)^i A^{n_1-1-i} [A^{n_2} B'_x \gamma B^* +  B \gamma (B^*)'_x (-A^*)^{n_2}] (A^*)^i  = \\
= \sum\limits_{i=0}^{n_1} (-1)^i A^{n_1-i} [ - A^{n_2} B \gamma\sigma_1^{-1}\sigma_2 B^* + B \sigma_2 \sigma_1^{-1}\gamma B^* (-A^*)^{n_2}]  (A^*)^i + \\
\quad + \sum\limits_{i=0}^{n_1-1} (-1)^i A^{n_1-1-i} [-A^{n_2+1} B \sigma_2\sigma_1^{-1}\gamma B^* -  B \gamma \sigma_1^{-1}\sigma_2 B^* A^* (-A^*)^{n_2}] (A^*)^i  = \\
= \text{using } \sigma_2\sigma_1^{-1}\gamma + \gamma \sigma_1^{-1}\sigma_2 = \sqrt{-1}\sigma_1 \\
= - \sqrt{-1} \sum\limits_{i=0}^{n_1} (-1)^i A^{n_1-i} A^{n_2} B \sigma_1 B^*(A^*)^i + (-1)^{n_1} A^{n_2} B \sigma_2\sigma_1^{-1}\gamma B^* (A^*)^{n_1} + \\
\quad + \sqrt{-1} \sum\limits_{i=0}^{n_1} (-1)^i A^{n_1-i} B\sigma_1 B^* (-A^*)^{n_2}  (A^*)^i - (-1)^{n_2} A^{n_1} B \gamma \sigma_1^{-1}\sigma_2 B^* (A^*)^{n_2}.
\end{array} \]
Finally, starting from the Lyapunov equation \eqref{eq:Lyapunov} and applying to this summation $\sum\limits_{i=0}^{n} A^{n-i} ..... (A^*)^i$ we obtain
from
\[ \sum\limits_{i=0}^{n} (-1)^i A^{n-i} [A \mathbb X + \mathbb X A^* + B\sigma_1B^* ] (A^*)^i = 0,
\]
that
\[  \sum\limits_{i=0}^{n} (-1)^i A^{n-i+1} \mathbb X (A^*)^i + \sum\limits_{i=0}^{n} (-1)^i A^{n-i} \mathbb X(A^*)^{i+1} =
- \sum\limits_{i=0}^{n} (-1)^i A^{n-i}  B\sigma_1B^* (A^*)^i.
\]
Since the first two sums are telescopic, we obtain
\begin{equation} \label{eq:Convolution}
  A^{n+1} \mathbb X + (-1)^{n} \mathbb X (A^*)^{n+1} =
- \sum\limits_{i=0}^{n} (-1)^i A^{n-i}  B\sigma_1B^* (A^*)^i.
\end{equation}
Plugging this formula into the last expression, we obtain
\[ \begin{array}{llll}
- \sqrt{-1} \sum\limits_{i=0}^{n_1} (-1)^i A^{n_1-i} A^{n_2} B \sigma_1 B^*(A^*)^i + (-1)^{n_1} A^{n_2} B \sigma_2\sigma_1^{-1}\gamma B^* (A^*)^{n_1} + \\
\quad + \sqrt{-1} \sum\limits_{i=0}^{n_1} (-1)^i A^{n_1-i} B\sigma_1 B^* (-A^*)^{n_2} (A^*)^i - (-1)^{n_2} A^{n_1} B \gamma \sigma_1^{-1}\sigma_2 B^* (A^*)^{n_2} = \\
=  \sqrt{-1} A^{n_2} [ A^{n_1+1} \mathbb X + (-1)^{n_1} \mathbb X (A^*)^{n_1+1}]  + (-1)^{n_1} A^{n_2} B \sigma_2\sigma_1^{-1}\gamma B^* (A^*)^{n_1} - \\
\quad - \sqrt{-1}[ A^{n_1+1} \mathbb X + (-1)^{n_1} \mathbb X (A^*)^{n_1+1}]  (-A^*)^{n_2}  - (-1)^{n_2} A^{n_1} B \gamma \sigma_1^{-1}\sigma_2 B^* (A^*)^{n_2} = \\
= \sqrt{-1} A^{n_2+n_1+1} \mathbb X  + \sqrt{-1} (-1)^{n_1}A^{n_2} \mathbb X (A^*)^{n_1+1}+ (-1)^{n_1} A^{n_2} B \sigma_2\sigma_1^{-1}\gamma B^* (A^*)^{n_1} - \\
 - \sqrt{-1} A^{n_1+1} \mathbb X (-A^*)^{n_2} - \sqrt{-1}(-1)^{n_1+n_2} \mathbb X (A^*)^{n_1+n_2} - (-1)^{n_2} A^{n_1} B \gamma \sigma_1^{-1}\sigma_2 B^* (A^*)^{n_2}.
\end{array} \]
Finally, exchanging here $n_1\leftrightarrow n_2$ and subtracting these two formulas we will obtain the following expression (notice that the terms depending
on $n_1+n_2$ are canceled)
\[ \begin{array}{llll}
\sqrt{-1} (-1)^{n_1}A^{n_2} \mathbb X (A^*)^{n_1+1} + (-1)^{n_1} A^{n_2} B \sigma_2\sigma_1^{-1}\gamma B^* (A^*)^{n_1} - \\
\quad \quad - \sqrt{-1} A^{n_1+1} \mathbb X (-A^*)^{n_2} - (-1)^{n_2} A^{n_1} B \gamma \sigma_1^{-1}\sigma_2 B^* (A^*)^{n_2} - \\
- \big[\sqrt{-1} (-1)^{n_2}A^{n_1} \mathbb X (A^*)^{n_2+1} + (-1)^{n_2} A^{n_1} B \sigma_2\sigma_1^{-1}\gamma B^* (A^*)^{n_2} - \\
\quad \quad - \sqrt{-1} A^{n_2+1} \mathbb X (-A^*)^{n_1} - (-1)^{n_1} A^{n_2} B \gamma \sigma_1^{-1}\sigma_2 B^* (A^*)^{n_1}\big] = \\
= - \sqrt{-1} (-1)^{n_2} A^{n_1} [ A\mathbb X - \sqrt{-1}B\gamma\sigma_1^{-1}\sigma_2 B^* + \mathbb XA^*-\sqrt{-1} B \sigma_2\sigma_1^{-1}\gamma B^*] (A^*)^{n_2} - \\
- \sqrt{-1} (-1)^{n_1} A^{n_2} [ A\mathbb X - \sqrt{-1}B\gamma\sigma_1^{-1}\sigma_2 B^* + \mathbb XA^*-\sqrt{-1} B \sigma_2\sigma_1^{-1}\gamma B^*] (A^*)^{n_1} = \\
- \sqrt{-1} (-1)^{n_2} A^{n_1} [ A\mathbb X + \mathbb XA^* + B \sigma_1 B^* ] (A^*)^{n_2} - \\
- \sqrt{-1} (-1)^{n_1} A^{n_2} [ A\mathbb X + \mathbb XA^* + B \sigma_1 B^* ] (A^*)^{n_1} = 0-0=0.
\end{array} \]
\qed

Under these conditions we will see that the evolutionary equation for $S$ and for $\beta$ will correspond to the KdV hierarchy.
Let us define a sequence $K_n$ of matrix-valued functions as follows:
\begin{defn} $K_n(x,t)$ is the $n$-th moments of the following function
\[ K(\lambda,x,t) = - S'_x(\lambda,x,t) S^{-1}(\lambda,x,t)\sigma_1^{-1} = \sum\limits_{i=0}^\infty \dfrac{K_n(x,t)}{\lambda^{n+1}}.
\]
\end{defn}
From this definition we are able to find a recursive formula for $K_n$'s involving the moments $H_n$ of $S(\lambda,x,t)$. Plugging the formula
\eqref{eq:DefHn} into Definition of $K_n$ we obtain
\[ \sum\limits_{i=0}^\infty \dfrac{K_n}{\lambda^{n+1}} \sigma_1 (I-\sum\limits_{i=0}^\infty \dfrac{H_n\sigma_1}{\lambda^{n+1}}) =
\sum\limits_{i=0}^\infty \dfrac{H'_i\sigma_1}{\lambda^{i+1}} 
\]
from where we get
\begin{eqnarray}
\label{eq:K0Def} K_0=H_0', \\
\label{eq:K1Def} K_1=H_1'+H_0'\sigma_1H_0, \\
\label{eq:KnRecursion} K_n = H_n' + \sum\limits_{i=0}^{n-1} K_i\sigma_1H_{n-1-i}.
\end{eqnarray}
\begin{lemma} The moments $K_n$ satisfy the following equality $K_n^* = (\sqrt{-1})^n K_n$.
\end{lemma}
\noindent\textbf{Proof:} Plugging $S^{-1*}(-\bar\lambda)$ from \eqref{eq:SSym} we find that
\[ \begin{array}{llll}
K^*(-\bar\lambda) & = \sigma_1^{-1} S^{-1*}(-\bar\lambda) (S^{-1*}(-\bar\lambda)\sigma_1 S^{-1}(-\bar\lambda))'_x \\
& = \sigma_1^{-1} (\sigma_1 S(\lambda)\sigma_1^{-1}) (\sigma_1 S(\lambda)\sigma_1^{-1} )'_x \\
& = - S(\lambda) S^{-1} (\lambda) (S(\lambda))'_x  S^{-1} (\lambda) \sigma_1^{-1} \\
& = -K(\lambda).
\end{array}\]
and the result follows. \qed

We will see evantually that the moments $K_n$ implement the family of solutions \eqref{eq:BtKdVn} by choosing $P_0 \sigma_1 = (\sqrt{-1})^n K_n$ in the formula
\eqref{eq:EvolKdVHirhy}. As a result for the understanding of the evolution equation \eqref{eq:EvolKdVHirhy}, we need
to understand the entries of 
\[ K_n' - \sigma_1^{-1} \gamma_* K_n + K_n \gamma_*\sigma_1^{-1} \quad ( \text{to be equal to } \sigma_1^{-1} (\gamma_*)'_t \sigma_1^{-1}).
\]
Notice that using $S'_x$ appearing in \eqref{eq:DS}
\[ K = - S'_x S^{-1}\sigma_1^{-1} = -\sigma_1^{-1}(\sigma_2\lambda+\gamma_*)\sigma_1^{-1} + S \sigma_1^{-1}(\sigma_2\lambda+\gamma)S^{-1}\sigma_1^{-1}.
\]
and plugging this formula for $K$ into $K_n' - \sigma_1^{-1} \gamma_* K_n + K_n \gamma_*\sigma_1^{-1}$ we find, denoting
$\Gamma=\sigma_1^{-1}(\lambda\sigma_2+\gamma)$, that
\begin{multline*} K_n' - \sigma_1^{-1} \gamma_* K_n + K_n \gamma_*\sigma_1^{-1} = \\
= -\sigma_1^{-1} (\gamma_*)'_x\sigma_1^{-1} + \lambda \sigma_1^{-1}\sigma_2 S\Gamma S^{-1}\sigma_1^{-1} - \lambda S\Gamma S^{-1} \sigma_1^{-1}\sigma_2 \sigma_1^{-1} +
\lambda \sigma_1^{-1} \gamma_* \sigma_1^{-1}\sigma_2 \sigma_1^{-1} - \lambda \sigma_1^{-1}\sigma_2 \sigma_1^{-1} \gamma_* \sigma_1^{-1}.
\end{multline*}
And thus it is necessary to understand the term $S\Gamma S^{-1}$, since all the other coefficients are in $\Rstar$. Simple calculations show that
\[ [S\Gamma S^{-1}]'_x = \Gamma_* S\Gamma S^{-1} - S\Gamma S^{-1}\Gamma_*, \quad \Gamma_* = \sigma_1^{-1}(\lambda \sigma_2 + \gamma_*).
\]
Denote the moments of $S\Gamma S^{-1}$ as $\bbmatrix{a_n&b_n\\c_n&d_n}$, and equate the coefficients of powers of $\lambda$ to find the following relations between
these moments:
\begin{multline*}
\bbmatrix{a_n'&b_n'\\c_n'&d_n'} = \sigma_1^{-1}\sigma_2 \bbmatrix {a_{n+1}&b_{n+1}\\c_{n+1}&d_{n+1}} -
\bbmatrix {a_{n+1}&b_{n+1}\\c_{n+1}&d_{n+1}} \sigma_1^{-1}\sigma_2 + \\
+ \sigma_1^{-1}\gamma_* \bbmatrix{a_n&b_n\\c_n&d_n} - \bbmatrix{a_n&b_n\\c_n&d_n} \sigma_1^{-1}\gamma_*
\end{multline*}
or the following recursive relations
\[ \left\{ \begin{array}{lllll}
a'_n & = - b_{n+1} + i c_n + i (\beta'-\beta^2) b_n, \\
b'_n & = 2 \beta b_n + i (d_n-a_n), \\
c'_n & = (a_{n+1} - d_{n+1}) + i (\beta'-\beta^2) (d_n-a_n) - 2 \beta c_n, \\
d'_n & =  b_{n+1} - i c_n - i (\beta'-\beta^2) b_n.
\end{array} \right. \]
Notice that from the first and the last equations it follows that $a'_{n+1}+d'_{n+1}=0$, and we will require
\begin{defn}(Normalization condition) $d_n = -a_n$.
\end{defn}
Notice that this must be checked for the initial value only. Under this condition, the last system of equations may be rewritten as follows
\[ \left\{ \begin{array}{lllll}
b_{n+1} & = - a'_n + i c_n + i (\beta'-\beta^2) b_n, \\
2 a_{n+1} & =  c'_n + 2 i (\beta'-\beta^2) a_n + 2 \beta c_n, \\
b'_n & = 2 \beta b_n - 2 i a_n.
\end{array} \right. \]
We can see here that actually the entries $a_{n+1}, b_{n+1}$ are recursively defined from the previous ones, but the last condition also imposes an
additional relation between them:
\begin{thm} Under normalization condition for $1,2,\ldots,n-1$, the entries of $b_n$'s satisfy the following recursive differential relation
\begin{equation} \label{eq:dbn+1}
4 b_{n+1}' = -i b_n''' + 4i (\beta'b_n)'.
\end{equation}
\end{thm}
\noindent\textbf{Remark:} This formula is very similar to the formula of deriving conservative functionals of the KdV equation
\cite[Formula 3.3]{bib:LaxAlmPeriod}, \cite{bib:GGKMVI}.

\noindent\textbf{Proof:} Let us rewrite the last equation in the system as
\[ 2 \beta b_{n+1} - 2 i a_{n+1} = b'_{n+1},
\]
where we plug the recursive formulas for $a_{n+1}, b_{n+1}$. The we obtain
\[ 2 \beta(- a'_n + i c_n + i (\beta'-\beta^2) b_n) - i ( c'_n + 2 i (\beta'-\beta^2) a_n + 2 \beta c_n) =
- a''_n + i c'_n + i [(\beta'-\beta^2) b_n]'.
\]
Rearranging the terms and performing some calculations this can be rewritten as
\[ 2ic_n' = a_n''-i[(\beta'-2\beta^2) b_n]' - i\beta b_n'' + i\beta' b_n'.
\]
Plugging here $2a_n = i(b_n'-2\beta b_n)$, obtained from the last equation in the system we obtain:
\[ 4ic_n' = ib_n'''+4i([\beta^2 b_n]' - \beta b_n''-\beta'b_n'-\beta''b_n).
\]
Finally, plugging this last formula for $4ic_n'$ in the derivative of the first formula of the system (recursive formula for $b_{n+1}$), we obtain
\[ \begin{array}{lll}
4b_{n+1}' & = -2i (b_n'-2\beta b_n)'' + ib_n'''+4i([\beta^2 b_n]' - \beta b_n''-\beta'b_n'-\beta''b_n) +
4i [(\beta'-\beta^2) b_n]' = \\
& = -i b_n''' + 4i (\beta'b_n)',
\end{array} \]
after cancellations. \qed

It is remained to show that the evolution, defined by \eqref{eq:BtKdVn} is implemented by $K_n$.

\begin{thm} \label{thm:nKdV}
For the evolutionary equation \eqref{eq:BtKdVn} of type $n$ the following formulas hold
\begin{eqnarray}
\label{eq:DHKdVHier} (H_N)'_t = (\sqrt{-1})^n [ (H_{N+n})'_x + K_0 \sigma_1 H_{N-1} + \ldots K_N H_N, \\
\label{eq:DStKdVHier} S'_t = (\sqrt{-1}\lambda)^n S'_x + \sqrt{-1}^n[\lambda^{n-1} K_0\sigma_1  + \lambda^{n-2} K_1 + \ldots + K_n] \sigma_1 S,
\end{eqnarray}
and as a result, the formula for the evolution of $\beta$ is
\begin{equation}
\label{eq:DbetatKdVHier} \sigma_1^{-1} (\gamma_*)'_t = i^n(K_n)'_x\sigma_1 + i^n [K_n\sigma_1,\sigma_1^{-1}\gamma_*],
\end{equation}
which can be solved recursively using \eqref{eq:dbn+1} starting from $b_0 = -\dfrac{1}{4} \beta'''_{xxx} + \dfrac{3}{2} (\beta'_x)^2$.
\end{thm}
\noindent\textbf{Proof:} Notice first that the formula \eqref{eq:DStKdVHier} follows immediately from \eqref{eq:DHKdVHier} at least for big powers of $\lambda$.
It is a simple check that it also holds for the small coefficients too.

So let us show that the formula \eqref{eq:DHKdVHier} holds. Using \eqref{eq:BtKdVn} and \eqref{eq:nDXt}
\[ \begin{array}{lllllll}
(H_N)'_t = (B^*\mathbb X^{-1}A^N B)'_t = (B^*)'_t \mathbb X^{-1} A^N B + B^*(\mathbb X^{-1})'_t A^N B + B^*\mathbb X^{-1}A^N B'_t = \\
 = (B^*)'_x (-\sqrt{-1}A^*)^n \mathbb X^{-1} A^N B + B^* \mathbb X^{-1} (\mathbb Y_n + \mathbb Y_{n-1}) \mathbb X^{-1} A^N B  +
B^*\mathbb X^{-1}A^N (\sqrt{-1}A)^n B'_x
\end{array} \]
Plugging the formula \eqref{eq:defY} for $\mathbb Y_{n-1},\mathbb Y_n$ and definitions for $m_n, m_{n-1}$, the last formula becomes
\begin{multline*} 
\dfrac{(H_N)'_t}{(\sqrt{-1})^n} =(B^*)'_x (-A^*)^n \mathbb X^{-1} A^N B - \sum\limits_{i=0}^n (-1)^i B^* \mathbb X^{-1} A^{n-i} B \sigma_2 B^* (A^*)^i \mathbb X^{-1} A^N B - \\
- \sum\limits_{i=0}^{n-1} (-1)^i B^* \mathbb X^{-1} A^{n-1-i} B \gamma B^* (A^*)^i \mathbb X^{-1} A^N B  +
B^*\mathbb X^{-1}A^N A^n B'_x.
\end{multline*}
From the formula \eqref{eq:Convolution} it follows that
\[ \mathbb X^{-1} A^n + (-1)^{n-1} (A^*)^n \mathbb X^{-1} = - \mathbb X^{-1}( \sum\limits_{i=0}^{n-1} A^{n-1-i}\sigma_1(A^*)^i)\mathbb X^{-1}
\]
Plugging this formula into the last expression we obtain:
\begin{multline*} 
\dfrac{(H_N)'_t}{(\sqrt{-1})^n} = (B^*)'_x(\mathbb X^{-1} A^N + \mathbb X^{-1}( \sum\limits_{i=0}^{n-1} A^{n-1-i}\sigma_1(A^*)^i)\mathbb X^{-1})A^N B - \\
 - \sum\limits_{i=0}^{n-1} B^*\mathbb X^{-1} A (AB\sigma_1+B\gamma)B^*(A^*)^i \mathbb X^{-1}A^NB -
 (-1)^n B^* \mathbb X^{-1} B\sigma_2B^*(A^*)^n\mathbb X^{-1} A^NB + B^*\mathbb X^{-1}A^{N+n} B'_x = \\
= (H_{N+n})'_x + \sum\limits_{i=0}^{n-1} (H_{n-1-i})'_x(-1)^i \sigma_1 B^* (A^*)^i\mathbb X^{-1}A^NB,
\end{multline*}
after some regrouping. So, it is remained to prove that
\[ (H_{N+n})'_x + \sum\limits_{i=0}^{n-1} (H_{n-1-i})'_x(-1)^i \sigma_1 B^* (A^*)^i\mathbb X^{-1}A^NB =
\sum\limits_{i=0}^n K_i \sigma_1 H_{N+n-i}.
\]
But this follows from the calculation of the coefficient of powers of $\lambda^{-n-1}\mu^{-N-1}$ at the expression
\[ S'_x(\lambda) S^{-1}(\lambda) S(\mu)
\]
in two ways. First from the definition of $K_n$ as the coefficient of $S'_x(\lambda) S^{-1}(\lambda) \sigma_1$ it follows that the coefficient of 
$\lambda^{-n-1}\mu^{-N-1}$ is
\[ \sum\limits_{i=0}^n K_i \sigma_1 H_{N+n-i}.
\]
On the other hand, Since $S^{-1}(\lambda) \sigma_1 = \sigma_1^{-1} S^*(-\bar\lambda)\sigma_1$ and since using \eqref{eq:Lyapunov} it is easy to show that
\[  S^*(-\bar\lambda)\sigma_1 S(\mu) - \sigma_1 = (-\lambda+\mu) \sigma_1 B^* (-\lambda I-A^*)^{-1} \mathbb X^{-1} (\mu I -A)^{-1} B\sigma_1,
\]
we will obtain the desired result. \qed
\subsection{\label{sec:Solitons}Solitons as examples of solutions}
One of the advantages of this approach is that we are able to create solutions of the $n$-th KdV equation.
For example taking a finite dimensional space $\mathcal H=\mathbb C^3$ and defining the operators as follows
\[ \begin{array}{lllll}
A & = \diag(-ik_j^2) = \bbmatrix{-ik_1^2&0&0\\0&-ik_2^2&0\\0&0&-ik_3^2} = -A^*, \\ 
B(x,t) & = \bbmatrix{e^{k_1x+k_1^{2n+1}t} b_1&0&0\\0&e^{k_2x+k_2^{2n+1}t} b_2&0\\0&0&e^{k_3x+k_3^{2n+1} t} b_3} \bbmatrix{1&ik_1\\1&ik_2\\1&ik_3},\\
\mathbb X(x,t) & = I + [\dfrac{e^{(k_i+k_j)x+(k_i^{2n+1}+k_j^{2n+1})t}}{k_i+k_j}b_ib^*_j] .
\end{array} \]
we will obtain that $B(x,t)$ satisfies differential equations \eqref{eq:DB}, \eqref{eq:nDBt} and $\mathbb X(x,t)$ satisfies \eqref{eq:DX}, \eqref{eq:nDXt}.
In other words, $\tau(x)=\det(\mathbb X(x,t))$ creates a soliton for the $n$-th KdV equation via formula \eqref{eq:BetaTau}.

Taking $\mathcal H$ in a different manner, we can obtain other solutions of \eqref{eq:nKdV}, similarly to the solution of the equation \eqref{eq:KdV}, presented in \cite{bib:KdVVessels}.

\bibliographystyle{alpha}
\bibliography{../../biblio}

\end{document}